\documentclass[amssymb,11pt]{amsart}
\usepackage{amscd,amssymb,euscript,amsthm}
\theoremstyle{plain}
\newtheorem{thm}{Theorem}[section] 
\newtheorem{cor}[thm]{Corollary}
\newtheorem{prop}[thm]{Proposition}

\theoremstyle{definition}

\theoremstyle{remark}
\newtheorem{rem}[thm]{Remark}

\numberwithin{equation}{section}
\def\<{\left<}
\def\>{\right>}
\def\cstar{$C^*$-algebra}
\begin{document}
\title{Dilation theory yesterday and today}
\author{William Arveson}
%
%
\address{Department of Mathematics,
University of California, Berkeley, CA 94720}
\email{arveson@math.berkeley.edu}
%
\date{22 February, 2009}

\begin{abstract} Paul Halmos' work in dilation theory began with 
a question and its answer: 
Which operators on a Hilbert space $H$ 
can be extended to normal operators 
on a larger Hilbert space $K\supseteq H$?  The answer is interesting and 
subtle.  

The idea of representing operator-theoretic 
structures in terms of conceptually simpler structures acting on larger Hilbert spaces has 
become a central one in the development of operator theory and, more 
generally, noncommutative analysis.  The work continues today.  In this article we summarize 
some of these diverse results and their history.  
\end{abstract}

\maketitle

\section{Preface}\label{S:in}

What follows is a brief account of the development of dilation theory that 
highlights Halmos' contribution to the circle of ideas.  
The treatment is not comprehensive.  I have chosen topics that have interested me over the years, 
while perhaps neglecting others.  
In order of 
appearance, the cast includes dilation theory 
for subnormal operators, operator valued measures and  
contractions, operator spaces, the role of extensions 
in dilation theory, commuting sets of operators, and 
semigroups of completely positive maps.  We emphasize Stinespring's theorem, 
but barely mention the model theory of 
Sz.-Nagy and Foias or its application to systems theory.    

After reflection on the common underpinnings of these 
results, it seemed a good idea to highlight 
the role of Banach $*$-algebras in their proofs, and I have done that.  
An appendix is included that summarizes what is needed.   
Finally, we have tried to make the subject accessible to students 
by keeping the prerequisites to a minimum; but of course we do assume familiarity 
with the basic theory of operators on Hilbert spaces and \cstar s.

\section{origins}\label{S:or}

Hilbert spaces are important because positive definite functions give rise to inner products
on vector spaces -- whose completions are Hilbert spaces --  
and positive definite functions are found in every corner of mathematics and mathematical physics.  
This association of a Hilbert space with a positive definite function involves 
a construction, and like all constructions that begin with objects in one category 
and generate objects in another category, it is best understood when viewed as a {\em functor}.  
We begin by discussing the properties of this functor in some detail since, while here they 
are simple and elementary,  similar properties will re-emerge later in other contexts.  

Let $X$ be a set and let 
$$
u: X\times X\to \mathbb C
$$
be 
a complex-valued function of two variables that is {\em positive definite} in the sense 
that for every $n=1,2,\dots$, 
every $x_1,\dots, x_n\in X$ and every set $\lambda_1, \dots, \lambda_n$ 
of complex numbers, one has 
\begin{equation}\label{orEq1}
\sum_{k,j=1}^n u(x_k,x_j)\bar\lambda_j\lambda_k\geq 0.  
\end{equation}
Notice that if $f: X\to H$ is a function from $X$ to a Hilbert space $H$  
with inner product $\langle\cdot,\cdot\rangle$, 
then the function $u: X\times X\to \mathbb C$ defined by  
\begin{equation}\label{orEq2}
u(x,y)=\langle f(x), f(y)\rangle, \qquad x,y\in X  
\end{equation}
is positive definite.  
By passing to a subspace 
of $H$ if necessary, one can obviously arrange that $H$ is the closed linear 
span of the set of vectors $f(X)$ in the range of $f$, and in that case the function $f: X\to H$ 
is said to be {\em minimal} (for the positive definite function $u$).  
Let us agree to say that two  
Hilbert space valued functions 
$f_1: X\to H_1$ and $f_2: X\to H_2$ are {\em isomorphic} if 
there is a unitary operator $U: H_1\to H_2$ such that 
$$
U(f_1(x))=f_2(x),\qquad x\in X.     
$$
A simple argument shows that all minimal functions for $u$ are isomorphic.  

For any positive definite function $u: X\times X\to \mathbb C$, 
a self-map $\phi: X\to X$ may or may not preserve the values of $u$ in the sense that 
$$
u(\phi(x),\phi(y))=u(x,y),\qquad x,y\in X; 
$$
but when this formula does hold, one would expect that  $\phi$ should 
acquire a Hilbert space identity.  
In order to discuss that, let us think of Hilbert spaces as the objects of 
a category whose morphisms are isometries; thus, a homomorphism from $H_1$ to 
$H_2$ is a linear isometry $U\in\mathcal B(H_1, H_2)$.  
Positive definite functions are also the objects of 
a category, in which a homomorphism from $u_1: X_1\times X_1\to \mathbb C$ to 
$u_2: X_2\times X_2\to \mathbb C$ is  a function $\phi: X_1\to X_2$ 
that preserves the positive structure in the sense that
\begin{equation}\label{orEq3}
u_2(\phi(x),\phi(y))=u_1(x,y),\qquad x,y\in X_1.  
\end{equation}

Given a positive definite function $u: X\times X\to \mathbb C$, 
one can construct a Hilbert space $H(u)$ and a function $f: X\to H(u)$ as follows.  
Consider the vector 
space $\mathbb CX$ of all complex valued functions 
$\xi: X\to \mathbb C$ with the property that $\xi(x)=0$ for all but a finite number of 
$x\in X$.  We can define a 
sesquilinear form $\langle \cdot,\cdot\rangle$ on $\mathbb CX$ by way of 
$$
\langle \xi,\eta\rangle=\sum_{x,y\in X}u(x,y)\xi(x)\bar\eta(y), \qquad \xi,\eta\in \mathbb CX,   
$$
and one finds that $\langle\cdot,\cdot\rangle$ is positive semidefinite because of 
the hypothesis on $u$.  
An application of the Schwarz inequality shows that 
the set 
$$
N=\{\xi\in\mathbb CX: \langle\xi,\xi\rangle=0\}
$$
is in fact a linear subspace of $\mathbb CX$, 
so this sesquilinear form can be promoted naturally to an inner product on 
the quotient $\mathbb CX/N$.  The completion of the inner product space $\mathbb CX/N$ is a 
Hilbert space $H(u)$, and  we can define the sought-after function $f: X\to H(u)$ as follows:
\begin{equation}\label{orEq4}
f(x)=\delta_x+N, \qquad x\in X, 
\end{equation}
where $\delta_x$ is the characteristic function of the singleton $\{x\}$.   
By construction, $u(x,y)=\langle f(x), f(y)\rangle$.  Note too that 
this function $f$ is {\em minimal} for $u$.  While there are many (mutually isomorphic) 
minimal functions for $u$,  we  
fix attention on the minimal function (\ref{orEq4}) that we have constructed.  

Given two positive definite functions $u_k: X_k\times X_k\to \mathbb C$, 
$k=1,2$, choose a homomorphism from $u_1$ to $u_2$, namely 
a function $\phi: X_1\to X_2$ that satisfies (\ref{orEq3}).   
Notice that while the two functions $f_k:X_k\to H(u_k)$
$$
f_1(x)=\delta_x+N_1, \quad f_2(y)=\delta_y+N_2, \qquad x\in X_1, \quad y\in X_2
$$
need not be injective, we do have the relations 
$$
\langle f_2(\phi(x)),f_2(\phi(y))\rangle_{H(u_2)} =
u_2(\phi(x),\phi(y))=u_1(x,y)=\langle f_1(x), f_1(y)\rangle_{H(u_1)}, 
$$
holding for all $x,y\in X_1$.  Since $H(u_1)$ is spanned by $f_1(X_1)$, 
a familiar and elementary argument (that we omit) 
shows that there is a unique linear isometry $U_\phi: H(u_1)\to H(u_2)$ such 
that 
\begin{equation}\label{orEq5}
U_\phi(f_1(x))=f_2(\phi(x)), \qquad x\in X_1.  
\end{equation}
At this point, it is straightforward to verify that the expected composition 
formulas $U_{\phi_1}U_{\phi_2}=U_{\phi_1\circ\phi_2}$ hold in general, and 
we conclude:   

\begin{prop}\label{orProp1} The construction (\ref{orEq4}) gives rise to a 
covariant functor $(u,\phi)\to (H(u),U_\phi)$ from the category of positive definite functions 
on sets to 
the category of complex Hilbert spaces.  
\end{prop}

It is significant that if $X$ is a topological space and $u:X\times X\to \mathbb C$ is 
a {\em continuous} positive definite function, 
then the associated map $f: X\to H(u)$ of (\ref{orEq4}) is also  continuous.  Indeed, 
this is immediate from (\ref{orEq2}): 
$$
\|f(x)-f(y)\|^2 =u(x,x)+u(y,y)-u(x,y)-u(y,x), \qquad x,y\in X.  
$$
The functorial nature of Proposition \ref{orProp1} pays immediate dividends:    

\begin{rem}\label{orRem1}[Automorphisms]
Every positive definite function 
$$
u: X\times X\to\mathbb C
$$
has an associated 
group of internal symmetries, namely the group $G_u$ of all bijections 
$\phi: X\to X$ that preserve $u$ in the sense that 
$$
u(\phi(x),\phi(y))=u(x,y), \qquad x,y\in X.  
$$ 
Notice that Proposition \ref{orProp1} implies that 
this group of symmetries has 
a natural unitary representation $U:G_u\to \mathcal B(H(u))$ associated 
with it.  Indeed, for every $\phi\in G_u$, the unitary operator $U_\phi\in\mathcal B(H(u))$ 
is defined uniquely by 
$$
U_\phi(f(x))=f(\phi(x)), \qquad x\in X.  
$$
The properties of this unitary representation of the automorphism group of $u$ often 
reflect important features of the environment that produced $u$.  
\end{rem}

\vskip0.1in
\noindent
{\bf Examples:}  There are many examples of positive definite functions;  
some of the more popular are reproducing kernels 
associated with domains in $\mathbb C^n$.  
Here is another example that 
is important for quantum physics and happens to be one of my favorites.  Let $Z$ be a (finite or 
infinite dimensional) Hilbert space and consider the positive definite 
function $u: Z\times Z\to \mathbb C$ defined by 
\begin{equation*}
u(z,w)=e^{\langle z,w\rangle}, \qquad z,w\in Z.  
\end{equation*}
We will write the Hilbert space $H(u)$ defined by the construction of Proposition \ref{orProp1} 
as $e^Z$, since it 
can be identified as the symmetric Fock space over the one-particle space $Z$.  
We will not make that identification here, but we do write the natural 
function (\ref{orEq4}) from $Z$ to $e^Z$ as $f(z)=e^z$, $z\in Z$.  

One finds that the automorphism group of Remark \ref{orRem1} is the 
full unitary group $\mathcal U(Z)$ of $Z$.  Hence 
the functorial nature of the preceding construction leads immediately to 
a (strongly continuous) 
unitary representation $\Gamma$ of the unitary group $\mathcal U(Z)$ on the 
Hilbert space  $e^Z$.  
In explicit terms, for $U\in\mathcal U(Z)$, $\Gamma(U)$ is the unique 
unitary operator of $e^Z$ that satisfies 
$$
\Gamma(U)(e^z)=e^{Uz}, \qquad U\in \mathcal U(Z), \quad z\in Z.  
$$
The map $\Gamma$ is called {\em second quantization} in the physics literature.  It 
has the property that for every one-parameter unitary group $\{U_t: t\in \mathbb R\}$ 
acting on $Z$, there is a corresponding one-parameter 
unitary group $\{\Gamma(U_t):t\in\mathbb R\}$ that 
acts on the ``first quantized" Hilbert space $e^Z$.  Equivalently, for every 
self adjoint operator $A$ on $Z$, there is a corresponding ``second quantized" 
self adjoint operator $d\Gamma(A)$ on $e^Z$ that is uniquely 
defined by the formula 
$$
e^{i t d\Gamma(A))}=\Gamma(e^{itA}), \qquad t\in \mathbb R,   
$$
as one sees by applying Stone's theorem which characterizes the generators 
of strongly continuous one-parameter unitary groups.  

Finally, one can exploit the functorial nature of this construction further to obtain 
a natural representation of the canonical commutation relations on $e^Z$, but we will 
not pursue that here.  

\section{Positive linear maps on commutative $*$-algebras}\label{S:pl}

The results of Sections \ref{S:ap} and \ref{S:pd} on subnormal operators, 
positive operator valued measures and the dilation theory of contractions 
can all be based on a single 
dilation theorem for positive linear maps of commutative Banach $*$-algebras.  
That commutative theorem has 
a direct commutative proof.  But since we require a more general noncommutative dilation 
theorem in Section \ref{S:cp} that contains it as a special case, we avoid 
repetition by merely stating the commutative result in this section.  What we 
want to emphasize here 
is the unexpected appearance of complete positivity even in this commutative context, and 
the functorial nature of dilation theorems of this kind.  

A Banach $*$-algebra is a Banach algebra $\mathcal A$ that is endowed with 
an isometric involution -- an antilinear mapping $a\mapsto a^*$ of 
$\mathcal A$ into itself that satisfies $a^{**}=a$, $(ab)^*=b^*a^*$ and $\|a^*\|=\|a\|$.  
In this section we will be concerned with Banach $*$-algebras that are 
{\em commutative}, and which have a multiplicative 
unit $\mathbf 1$ that satisfies $\|\mathbf 1\|=1$. 
The basic properties 
of Banach $*$-algebras and their connections with \cstar s are summarized in Appendix \ref{S:a1}.  

An operator-valued linear map $\phi: \mathcal A\to \mathcal B(H)$ of a Banach 
$*$-algebra is said to be {\em positive} if $\phi(a^*a)\geq 0$ for every $a\in \mathcal A$.  
The most important fact about operator-valued positive linear maps of commutative 
algebras is something of a miracle.  It asserts that a positive linear map 
$\phi:\mathcal A\to\mathcal B(H)$ of a commutative Banach $*$-algebra $\mathcal A$ is 
{\em completely} positive in the following sense: For every $n$-tuple  
$a_1, \dots, a_n$ of elements of $\mathcal A$, the $n\times n$ operator matrix 
$(\phi(a_i^*a_j))$ is positive in the natural sense that for every $n$-tuple of 
vectors $\xi_1,\dots, \xi_n\in H$,  one has 
\begin{equation}\label{plEq1}
\sum_{i,j=1}^n\langle \phi(a_i^*a_j)\xi_j,\xi_i\rangle \geq 0.  
\end{equation}
Notice that the hypothesis $\phi(a^*a)\geq 0$ is the content of 
these inequalities for the special case $n=1$.  
This result for commutative \cstar s $\mathcal A$ 
is due to Stinespring (see Theorem 4 of \cite{stine}), 
and the proof of (\ref{plEq1}) can be based on 
that result  combined with the 
properties of the completion map $\iota: \mathcal A\to C^*(\mathcal A)$ 
that carries a commutative Banach $*$-algebra $\mathcal A$ to its enveloping \cstar\ 
$C^*(\mathcal A)\cong C(X)$ (see Remark \ref{a1Rem1} of Appendix \ref{S:a1}).  

The notion of complete positivity properly belongs to the noncommutative world.  
We will return to it in Section \ref{S:cp} where we will prove a general 
result (Theorem \ref{cpThm1}) which, 
when combined with (\ref{plEq1}), implies the following assertion about 
positive linear maps of commutative $*$-algebras.    

\vskip0.1in
\noindent
{\bf Scholium A:}\label{plSc1}
{\em Let $\mathcal A$ be a commutative Banach $*$-algebra with unit and 
let $H$ be a Hilbert space.  For every 
operator-valued linear map $\phi: \mathcal A\to \mathcal B(H)$ satisfying 
$\phi(a^*a)\geq 0$, $a\in \mathcal A$, there is a pair $(V,\pi)$ consisting 
of a representation $\pi: \mathcal A\to\mathcal B(K)$ of $\mathcal A$ on another 
Hilbert space $K$ and a linear operator $V\in\mathcal B(H,K)$ such that 
\begin{equation}\label{plEq2}
\phi(a)=V^*\pi(a)V, \qquad a\in \mathcal A.     
\end{equation}
Moreover, $\phi$ is necessarily bounded, its norm is given by 
\begin{equation}\label{plSc3}
\sup_{\|a\|\leq 1}\|\phi(a)\|=\|\phi(\mathbf 1)\|=\|V\|^2,   
\end{equation}
and $V$ can be taken to be an isometry when $\phi(\mathbf 1)=\mathbf 1$.  
}

\begin{rem}[Minimality and uniqueness of dilation pairs]\label{plRem1}
Fix $\mathcal A$ as above.  By a {\em dilation pair} for $\mathcal A$ we mean a pair 
$(V,\pi)$ consisting of a representation $\pi: \mathcal A\to \mathcal B(K)$ 
and a bounded linear map $V: H\to K$ from some other Hilbert space $H$ into 
the space $K$ on which $\pi$ acts.  
A dilation pair $(V,\pi)$ is said to be {\em minimal} if 
the set of vectors $\{\pi(a)V\xi: a\in \mathcal A, \ \xi\in H\}$ has 
$K$ as its closed linear span.  By replacing $K$ with an appropriate 
subspace and $\pi$ with an appropriate subrepresentation, we can 
obviously replace every such pair with a minimal one.  Moreover, the 
representation associated with a minimal pair must be nondegenerate, 
and therefore $\pi(\mathbf 1)=\mathbf 1_K$.     

Note that 
every dilation pair $(V,\pi)$ gives rise 
to a positive linear map $\phi: \mathcal A\to \mathcal B(H)$ 
that is defined by the formula (\ref{plEq2}), 
and we say that $(V,\pi)$ is a dilation pair {\em for} $\phi$.  A positive map 
$\phi$ has many dilation pairs associated with it, but the minimal ones 
 are equivalent 
in the following sense: If $(V_1, \pi_1)$ and $(V_2, \pi_2)$ are two {\em minimal} dilation 
pairs for $\phi$ then there is a unique unitary operator $W: K_1\to K_2$ such that 
\begin{equation}\label{plEq4}
WV_1=V_2, \quad {\rm{and}}\quad W\pi_1(a)=\pi_2(a)W, \qquad a\in \mathcal A.  
\end{equation}
The proof amounts to little more than 
checking inner products on the two generating sets $\pi_1(\mathcal A)V_1H\subseteq K_1$ 
and $\pi_2(\mathcal A)V_2H\subseteq K_2$ and noting that 
\begin{align*}
\langle \pi_2(a)V_2\xi,\pi_2(b)V_2\eta\rangle&=
\langle \pi_2(b^*a)V_2\xi,\eta\rangle=\langle\phi(b^*a)\xi,\eta\rangle
\\
&=
\langle \pi_1(a)V_1\xi,\pi_1(b)V_1\eta\rangle,   
\end{align*}
for $a,b\in\mathcal A$ and $\xi,\eta\in H$.  

Finally, note that in cases where $\phi(\mathbf 1)=\mathbf 1$, the operator $V$ of 
a minimal pair $(V,\pi)$ is an isometry, so by making an obvious identification we 
can replace $(V,\pi)$ with an equivalent one in which $V$ is the inclusion map 
of $H$ into a larger Hilbert space $\iota: H\subseteq K$ and $\pi$ is a representation 
of $\mathcal A$ on $K$.  After these identifications, 
(\ref{plEq2}) reduces to the more traditional assertion 
\begin{equation}\label{plEq5}
\phi(a)=P_H\pi(a)\restriction_H, \qquad a\in \mathcal A.  
\end{equation}
\end{rem}

\begin{rem}[Functoriality]\label{cdRem1}  
It is a worthwhile exercise to think carefully about what a dilation 
actually {\em is}, and the way to do that is to think 
in categorical terms.   Fix a commutative Banach $*$-algebra $\mathcal A$ 
with unit $\mathbf 1$.  Operator-valued positive linear maps of $\mathcal A$ 
are the objects of a category, in which a homomorphism from 
$\phi_1: \mathcal A\to \mathcal B(H_1)$ 
to $\phi_2: \mathcal A\to \mathcal B(H_2)$ is defined as a unitary operator 
$U: H_1\to H_2$ satisfying $U\phi_1(a)=\phi_2(a)U$ for all $a\in \mathcal A$; 
equivalently, $U$ should implement a unitary equivalence of positive linear maps of $\mathcal A$.   
Thus the positive linear maps of $\mathcal A$ can be viewed as 
a groupoid -- a category in which every 
arrow is invertible.  

There is a corresponding groupoid whose objects are minimal dilation pairs $(V,\pi)$.  
Homomorphisms of dilation pairs $(V_1, \pi_1)\to (V_2, \pi_2)$ (here 
$\pi_j$ is a representation of $\mathcal A$ on $K_j$ and $V_j$ is an 
operator in $\mathcal B(H_j,K_j)$) are defined as unitary 
operators $W: K_1\to K_2$ that satisfy 
$$
W\pi_1(a)=\pi_2(a)W, \qquad a\in \mathcal A, \quad {\rm{and}}\quad WV_1=V_2.  
$$ 
The ``set" of all dilation pairs 
for a fixed positive linear map $\phi: \mathcal A\to \mathcal B(H)$ is a subgroupoid, 
and we have already seen in Remark \ref{plRem1} that its elements are all 
isomorphic.  But here we are mainly concerned with how the dilation functor 
treats arrows between different positive linear maps.  

A functor is the end product of a {\em construction}.  
In order to describe how the dilation functor acts on arrows, we need 
more information than the statement of Scholium A contains, namely 
the following: {\em There is a construction which starts with a positive linear map 
$\phi: \mathcal A\to \mathcal B(H)$ and generates a particular dilation 
pair $(V,\pi)_\phi$ from that data.}  Scholium A asserts that such 
dilation pairs exist for every $\phi$, but since the proof is missing,  we have not 
seen the construction.  
Later on, however, we will show how to construct a particular dilation 
pair $(V,\pi)_\phi$ from 
a completely positive map $\phi$ when we prove Stinespring's theorem in section 
\ref{S:cp}.  That construction 
is analogous to the construction underlying (\ref{orEq4}), which exhibits an 
explicit 
function $f: X\to H(u)$ that arises from the construction of 
the Hilbert space $H(u)$, starting with a positive definite function $u$.  
In order to continue the current discussion, we ask the reader 
to assume the result of the construction of Theorem \ref{cpThm1}, namely 
that we are somehow given a {\em particular} 
dilation pair $(V,\pi)_\phi$ for every positive linear map 
$\phi: \mathcal A\to \mathcal B(H)$.  

That puts us in position to describe how the dilation functor acts on arrows.  Given 
two positive linear maps $\phi_j: \mathcal A\to \mathcal B(H_j)$, 
$j=1,2$, let $U: H_1\to H_2$ be a unitary operator satisfying 
$U\phi_1(a)=\phi_2(a)U$ for $a\in \mathcal A$.  Let $(V_1, \pi_1)$ 
and $(V_2,\pi_2)$ be the dilation pairs that have been 
constructed from $\phi_1$ and $\phi_2$ respectively.  Notice that since 
$U^*\phi_2(a)U=\phi_1(a)$ for $a\in \mathcal A$, it follows 
that $(V_2U, \pi_2)$ is a second minimal dilation pair for $\phi_1$.  
By (\ref{plEq4}), there is a unique unitary operator 
$\tilde U: K_1\to K_2$ that satisfies 
$$
\tilde UV_1=V_2U, \quad {\rm{and}}\quad \tilde U\pi_1(a)=\pi_2(a)\tilde U, \qquad a\in \mathcal A.  
$$
One can now check that the association $\phi,U \to (V,\pi)_\phi,\tilde U$ defines 
a covariant functor from the groupoid of 
operator-valued positive linear maps of $\mathcal A$ to the 
groupoid of minimal dilation pairs for $\mathcal A$.  
\end{rem}

\section{Subnormality}\label{S:ap}

An operator $A$ on a Hilbert space $H$ is said to be {\em subnormal} if it 
can be extended to a normal operator on a larger Hilbert space.  More 
precisely, there should exist a 
normal operator $B$ acting on a Hilbert space $K\supseteq H$ that leaves 
$H$ invariant and restricts to $A$ on $H$.  
Halmos' paper \cite{haldil} introduced the concept, and 
grew out of his observation that 
a subnormal operator $A\in\mathcal B(H)$ must satisfy the following 
system of peculiar inequalities: 
\begin{equation}\label{apEq0}
\sum_{i,j=0}^n\langle A^i\xi_j,A^j\xi_i\rangle\geq 0,\qquad \forall\ \xi_0,\xi_1,\dots, \xi_n\in H, 
\quad n=0,1,2,\dots.  
\end{equation}
It is an instructive exercise with inequalities involving $2\times 2$ operator matrices 
to show that the case $n=1$ of 
(\ref{apEq0}) is equivalent to the single 
operator inequality $A^*A\geq AA^*$, a property 
called {\em hyponormality} today.  Subnormal operators are certainly hyponormal, but 
the converse is false even for weighted shifts (see Problem 160 of \cite{halProbBook}).  
Halmos showed that the full set of inequalities (\ref{apEq0}) -- together with 
a second system of necessary inequalities that we do not reproduce here -- implies that $A$ 
is subnormal.  Several years later, his student J. Bram proved that the second system 
of inequalities follows from the first \cite{bram}, and 
simpler proofs of that fact based on semigroup considerations 
emerged later \cite{szfran}.   Hence the system of 
inequalities (\ref{apEq0}) is by itself necessary and sufficient for subnormality.  

It is not hard to reformulate Halmos' notion of subnormality (for single operators) 
in a more general way that applies to several operators.  
Let $\Sigma$ be a commutative semigroup (written additively) that contains a neutral element $0$.  
By a {\em representation} of $\Sigma$ we mean 
an operator valued function $s\in \Sigma\mapsto A(s)\in \mathcal B(H)$ 
satisfying $A(s+t)=A(s)A(t)$ and $A(0)=\mathbf 1$.  Notice that we make no assumption 
on the norms $\|A(s)\|$ as $s$ varies over $\Sigma$.  
For example, a commuting set $A_1,\dots, A_d$ of operators on a Hilbert space 
$H$ gives rise to a representation of the $d$-dimensional additive semigroup 
$$
\Sigma=\{(n_1,\dots, n_d): n_1\geq 0,\dots, n_d\geq 0\}
$$ 
by way of 
$$
A(n_1,\dots, n_d)=A_1^{n_1}\cdots A_d^{n_d}, 
\qquad (n_1,\dots, n_d)\in \Sigma. 
$$
In general, a representation $A:\Sigma\to \mathcal B(H)$ 
is said to be {\em subnormal} if there is a Hilbert space $K\supseteq H$ and a  
representation $B: \Sigma\to \mathcal B(K)$ consisting of normal operators 
such that each $B(s)$ leaves $H$ invariant and 
$$
B(s)\restriction_H=A(s), \qquad s\in \Sigma.
$$ 
We now apply Scholium A to prove a general statement 
about commutative operator semigroups that 
contains the Halmos-Bram characterization of subnormal operators, as well as higher dimensional 
variations of it that apply to semigroups generated by a finite or 
even infinite number of mutually commuting operators.  

\begin{thm}\label{apThm1} Let $\Sigma$ be a commutative semigroup with $0$.  
A representation $A: \Sigma\to \mathcal B(H)$ is subnormal iff for 
every $n\geq 1$, every $s_1,\dots, s_n\in \Sigma$ and every $\xi_1,\dots, \xi_n\in H$, 
one has 
\begin{equation}\label{apEq1}
\sum_{i,j=1}^n\langle A(s_i)\xi_j,A(s_j)\xi_i\rangle \geq 0.  
\end{equation}
\end{thm}

\begin{proof}  The proof that the system 
of inequalities (\ref{apEq1}) is necessary for subnormality is straightforward, and 
we omit it.  
Here we outline a proof of the converse, 
describing all essential steps in the construction but 
leaving routine calculations for the reader.  We shall make use of the 
hypothesis (\ref{apEq1}) in the following form: 
For every function $s\in\Sigma\mapsto \xi(s)\in H$ such that $\xi(s)$ vanishes for all but a finite 
number of $s\in \Sigma$,  one has 
\begin{equation}\label{apEq2}
\sum_{s,t\in \Sigma}\langle A(s)\xi(t), A(t)\xi(s)\rangle \geq 0.  
\end{equation}

We first construct an appropriate commutative Banach $*$-algebra.  
Note that the direct 
sum of semigroups $\Sigma\oplus \Sigma$ is a commutative semigroup with 
zero element $(0,0)$, but unlike $\Sigma$ it has a natural involution $x\mapsto x^*$ 
defined by $(s,t)^*=(t,s)$, $s,t\in \Sigma$.  We will also make use of a weight 
function $w: \Sigma\oplus\Sigma\to [1,\infty)$ defined as follows: 
$$
w(s,t)=\max(\|A(s)\|\cdot\|A(t)\|, 1), \qquad s,t\in \Sigma.   
$$
Straightforward verification (using $\|A(s+t)\|\leq \|A(s)\|\cdot\|A(t)\|$) 
establishes the properties 
$$
1\leq w(x+y)\leq w(x)w(y), \quad w(x^*)=w(x), \qquad x,y\in\Sigma\oplus\Sigma.  
$$ 
Note too that $w((0,0))=1$ because $A(0)=\mathbf 1$.  
Consider the Banach space $\mathcal A$ 
of all functions $f:\Sigma\oplus\Sigma\to \mathbb C$ 
having finite weighted $\ell^1$-norm
\begin{equation}\label{apEq2.1}
\|f\|=\sum_{x\in \Sigma\oplus\Sigma}|f(x)|\cdot w(x)<\infty.  
\end{equation}
Since $w\geq 1$, the norm on $\mathcal A$ dominates the 
ordinary $\ell^1$ norm, so that every function in $\mathcal A$ 
belongs to $\ell^1(\Sigma\oplus\Sigma)$.  
Ordinary convolution of functions defined on commutative semigroups 
$$
(f*g)(z)=\sum_{\{x,y\in\Sigma\oplus\Sigma: \ x+y=z\}}f(x)g(y), \qquad z\in \Sigma\oplus\Sigma
$$
defines an associative commutative multiplication in $\ell^1(\Sigma\oplus\Sigma)$, 
and it is easy to check that the 
above properties of the weight function $w$ imply that with 
respect to convolution and the involution $f^*(s,t)=\bar f(t,s)$,  $\mathcal A$ 
becomes a commutative 
Banach $*$-algebra with normalized unit $\delta_{(0,0)}$.  

We now use the semigroup $A(\cdot)$ to construct 
a linear map $\phi: \mathcal A\to \mathcal B(H)$: 
$$
\phi(f)=\sum_{(s,t)\in \Sigma\oplus\Sigma}f(s, t)A(s)^*A(t).     
$$
Note that $\|\phi(f)\|\leq \|f\|$ because of the definition  
norm of $f$ in terms of the weight function $w$. 
Obviously $\phi(\delta_{(s,t)})=A(s)^*A(t)$ for all $s,t\in \Sigma$, and 
in particular $\phi(\delta_{(0,0)})=\mathbf 1$.     
It is also obvious that $\phi(f^*)=\phi(f)^*$ for $f\in\mathcal A$.  

What is most important for us is that 
$\phi$ is a {\em positive} linear map, namely for every $f\in\mathcal A$ 
and every vector $\xi\in H$ 
\begin{equation}\label{apEq3}
 \langle\phi((f^*)*f)\xi,\xi\rangle \geq 0.  
\end{equation}
To deduce this from (\ref{apEq2}), note that 
since $\phi:\mathcal A\to\mathcal B(H)$ is a bounded linear map and every function in 
$\mathcal A$ can be 
norm-approximated by functions which are finitely nonzero, it suffices to 
verify (\ref{apEq3}) for functions $f:\Sigma\oplus\Sigma\to \mathbb C$ 
such that $f(x)=0$ for all but a finite number of $x\in \Sigma\oplus\Sigma$.  
But for two finitely supported functions 
$f,g\in \mathcal A$ and any function 
$H: \Sigma\oplus\Sigma\to \mathbb C$, the definition of convolution implies 
that $f*g$ is finitely supported, and 
$$
\sum_{z\in \Sigma\oplus\Sigma}(f*g)(z)H(z)=\sum_{x,y\in \Sigma\oplus\Sigma}f(x)g(y)H(x+y).  
$$
Fixing $\xi\in H$ and taking $H(s,t)=\langle A(s)^*A(t)\xi,\xi\rangle
=\langle A(t)\xi,A(s)\xi\rangle$, we conclude from the preceding formula that 
\begin{align*}
\langle \phi(f*g)\xi,\xi\rangle &=\sum_{s,t,u,v\in \Sigma}f(s,t)g(u,v)
\langle A(t+v)\xi,A(s+u)\xi\rangle 
\\
&=\sum_{s,t,u,v\in \Sigma}f(s,t)g(u,v)
\langle A(t)A(v)\xi,A(u)A(s)\xi\rangle .  
\end{align*}
Thus we can  write 
\begin{align*}
\langle \phi(f^**f)\xi,\xi\rangle &=
\sum_{s,t,u,v\in \Sigma}\bar f(t,s)f(u,v)
\langle A(t)A(v)\xi,A(u)A(s)\xi\rangle 
\\
&=
\sum_{t,u\in \Sigma}
\langle A(t)(\sum_{v\in\Sigma}f(u,v)A(v)\xi),
A(u)(\sum_{s\in\Sigma}f(t,s)A(s)\xi)\rangle 
\\
&=
\sum_{t,u\in\Sigma}\langle A(t)\xi(u),A(u)\xi(t)\rangle, 
\end{align*}
where $t\in\Sigma\mapsto \xi(t)\in H$ is the vector function 
$$
\xi(t)=\sum_{s\in \Sigma}f(t,s)A(s)\xi, \qquad t\in \Sigma.  
$$
Notice that the rearrangements of summations carried out in the preceding 
formula are legitimate because all sums are finite, and in fact 
the vector function $t\mapsto\xi(t)$ is itself finitely nonzero.  (\ref{apEq3}) now follows 
from (\ref{apEq2}).  
 
At this point, we can apply Scholium A to find a Hilbert space 
$K$ which contains $H$ and a {\em $*$-representation} $\pi: \mathcal A\to \mathcal B(K)$ 
such that 
$$
P_H\pi(f)\restriction_H =\phi(f),\qquad f\in \mathcal A.
$$
Hence the map $x\mapsto \pi(\delta_x)$ is a $*$-preserving representation  
of the $*$-semigroup $\Sigma\oplus\Sigma$, which can be further decomposed by way of 
$\pi(\delta_{(s,t)})=B(s)^*B(t)$, where $B:\Sigma\to\mathcal B(K)$ is  the  
representation $B(t)=\pi(\delta(0,t))$.  Since the 
commutative semigroup of operators $\{\pi(\delta_x): x\in \Sigma\oplus \Sigma\}$ is 
closed under the $*$-operation, 
$B(\Sigma)$ is a semigroup of mutually commuting normal operators.  After 
taking $s=0$ in 
the formulas $P_HB(s)^*B(t)\restriction_H=A(s)^*A(t)$, one finds that 
$A(t)$ is the compression of $B(t)$ to $H$.  Moreover, since for every $t\in\Sigma$
$$
P_HB(t)^*B(t)\restriction_H=A(t)^*A(t)=P_HB(t)^*P_HB(t)\restriction_H, 
$$
we have $P_HB(t)^*(\mathbf 1-P_H)B(t)P_H=0$.  Thus we have shown that $H$ is invariant 
under $B(t)$ and the restriction of $B(t)$ to $H$ is $A(t)$.   
\end{proof}

\begin{rem}[Minimality and functoriality]\label{plRem2}  Let $\Sigma$ 
be a commutative semigroup with zero.  
A normal extension $s\in\Sigma\mapsto B(s)\in \mathcal B(K)$ of 
a representation $s\in \Sigma\mapsto A(s)\in\mathcal B(H)$ on a Hilbert space $K\supseteq H$ is 
said to be {\em minimal} if the set of vectors 
$\{B(t)^*\xi: t\in \Sigma,\  \xi\in H\}$ has $K$ as its closed linear span.  
This corresponds to the notion of minimality described in Section \ref{S:cp}.  
The considerations of Remark 
\ref{plRem1} imply that all minimal dilations are equivalent,  and we can speak 
unambiguously of {\em the} minimal normal extension of $A$.  A similar comment applies to the 
functorial nature of the map which carries subnormal representations of 
$\Sigma$ to their minimal normal extensions.  
\end{rem}

\begin{rem}[Norms and flexibility]\label{apRem3}  It is a fact that the minimal 
normal extension $B$ of $A$ satisfies $\|B(t)\|=\|A(t)\|$ for $t\in \Sigma$.  The 
inequality $\geq$ is obvious since $A(t)$ is the restriction of $B(t)$ to an invariant subspace.  
However, if one attempts to use  
the obvious norm estimate for representations of Banach $*$-algebras (see Appendix \ref{S:a1} for 
more detail) to 
establish the opposite inequality, one finds that the above construction gives only 
$$
\|B(t)\|=\|\pi(\delta_{(0,t)})\|\leq \|\delta_{(0,t)}\|=w(0,t)=\max(\|A(t)\|, 1), 
$$
which is not good enough when $\|A(t)\|<1$.  On the other hand, we can use the flexibility 
in the possible norms of $\mathcal A$ to obtain the correct estimate 
as follows.  For each $\epsilon>0$, define a new weight function $w_\epsilon$ 
on $\Sigma\oplus\Sigma$ by 
$$
w_\epsilon(s,t)=\max(\|A(s)\|\cdot \|A(t)\|, \epsilon), \qquad s,t\in \Sigma.  
$$
If one uses $w_\epsilon$ in place of $w$ in the definition (\ref{apEq2.1}) of 
the norm on $\mathcal A$, 
one obtains another commutative Banach $*$-algebra which serves equally well as the original 
to construct the minimal normal dilation $B$ of $A$, and it has the additional property 
that $\|B(t)\|\leq \max(\|A(t)\|, \epsilon)$ for $t\in \Sigma$.  Since $\epsilon$ can 
be arbitrarily small, the desired estimate $\|B(t)\|\leq \|A(t)\|$ follows.  In particular, 
for every $t\in\Sigma$ we have $A(t)=0\implies B(t)=0$.  
\end{rem}

\section{Commutative dilation theory}\label{S:pd}

Dilation theory began with two papers of Naimark, written and published 
somehow during 
the darkest period of world war II: \cite{nai2}, \cite{nai1}.   Naimark's theorem asserts 
that a countably additive measure $E: \mathcal F\to \mathcal B(H)$ defined on a $\sigma$-algebra 
$\mathcal F$ of subsets of a set $X$ that takes 
values in the set of positive operators on a Hilbert space $H$ and satisfies 
$E(X)=\mathbf 1$ can be expressed 
in the form 
$$
E(S)=P_HQ(S)\restriction_H, \qquad S\in\mathcal F, 
$$
where $K$ is a Hilbert space containing $H$ 
and $Q: \mathcal F\to \mathcal B(K)$ is a spectral measure.  A version of 
Naimark's theorem (for regular Borel measures on topological spaces) can be 
found on p. 50 of \cite{paulsenBk2}.  Positive operator valued measures $E$ have 
become fashionable in quantum physics and quantum information theory, where 
they go by the  unpronounceable acronym POVM.  It is interesting that the Wikipedia page for  
projection-operator-valued-measures 
({\tt{http://en.wikipedia.org/wiki/POVM}}) 
contains more information about 
Naimark's famous theorem than the Wikipedia page for Naimark himself 
({\tt{http://en.wikipedia.org/wiki/Mark\_Naimark}}).  

In his subnormality paper \cite{haldil}, 
Halmos showed that every contraction $A\in \mathcal B(H)$ has a unitary dilation 
in the sense that there is a unitary operator $U$ acting on a larger Hilbert space 
$K\supseteq H$ such that 
$$
A = P_HU\restriction_H.  
$$  
Sz.-Nagy extended that in a most significant way \cite{szNdil} by showing that every 
contraction has a unitary {\em power} dilation, and the latter result 
ultimately became the cornerstone for 
an effective model theory for Hilbert space contractions \cite{nagFoias}.  Today, these results 
belong to the toolkit of every operator theorist, and can be found in many 
textbooks.  In this section we merely state Sz.-Nagy's theorem and sketch 
a proof that is in the spirit of the preceding discussion.  

\begin{thm}\label{pdThm1}  Let $A\in\mathcal B(H)$ be an operator satisfying 
$\|A\|\leq 1$.  Then there is a unitary operator $U$ acting on         a Hilbert space 
$K$ containing $H$ such that 
\begin{equation}\label{pdEq1}
A^n=P_HU^n\restriction_H, \qquad n=0,1,2,\dots.  
\end{equation}
If $U$ is {\em minimal} in the sense that $K$ is the closed linear span of 
$\cup_{n\in\mathbb Z}U^nH$, then it is uniquely determined up to a natural 
unitary equivalence.  
\end{thm}

\begin{proof}[Sketch of proof]  Consider the commutative Banach $*$-algebra 
$\mathcal A=\ell^1(\mathbb Z)$, with multiplication and involution given 
by 
$$
(f*g)(n)=\sum_{k=-\infty}^{+\infty} f(k)g(n-k), \quad f^*(n)=\bar f(-n), \qquad n\in\mathbb Z, 
$$
and normalized unit $\mathbf 1=\delta_0$.  Define $A(n)=A^n$ if $n\geq 0$ and 
$A(n)=A^{*|n|}$ if $n<0$.  Since $\|A(n)\|\leq 1$ for every $n$, we can define 
a linear map $\phi: \mathcal A\to \mathcal B(H)$ in the obvious way 
$$
\phi(f)=\sum_{n=-\infty}^{+\infty} f(n)A(n), \qquad f\in \mathcal A.  
$$
It is obvious that $\|\phi(f)\|\leq \|f\|$, $f\in \mathcal A$, but not at all obvious 
that $\phi$ is a positive linear map.  However, there is a standard method for showing that 
for every $\xi\in H$, the sequence of complex numbers $a_n=\langle A(n)\xi,\xi\rangle$, 
$n\in \mathbb Z$, is of positive type in the sense that for every finitely nonzero 
sequence of complex numbers $\lambda_n$, $n\in \mathbb Z$, one has 
$
\sum_{n\in \mathbb Z}a_{n-m}\lambda_n\bar\lambda_m\geq 0; 
$
for example, see p. 36 of \cite{paulsenBk2}.  By approximating $f\in \mathcal A$ 
in the norm of $\mathcal A$ 
with finitely nonzero functions and using 
$$
\langle\phi((f^*)*f)\xi,\xi\rangle =
\sum_{m,n=-\infty}^{+\infty}\langle A(n-m)\xi,\xi\rangle f(n)\bar f(m), 
$$ 
it follows that $\langle\phi(f)\xi,\xi\rangle\geq 0$, 
and we may conclude that $\phi$ is a positive linear map 
of $\mathcal A$ to $\mathcal B(H)$ satisfying 
$\phi(\delta_0)=\mathbf 1$.  

Scholium A implies that there is a $*$-representation $\pi: \mathcal A\to \mathcal B(K)$ 
of $\mathcal A$ on a larger Hilbert space $K$ such that $\pi(f)$ compresses to $\phi(f)$ 
for $f\in\mathcal A$.  Finally, 
since the enveloping \cstar\ of $\mathcal A=\ell^1(\mathbb Z)$ is the commutative 
\cstar\ $C(\mathbb T)$, the representation $\pi$ promotes to a representation 
 $\tilde\pi: C(\mathbb T)\to\mathcal B(K)$ (see Appendix \ref{S:a1}).  
 Taking $z\in C(\mathbb T)$ to 
be the coordinate variable, we obtain a unitary operator $U\in \mathcal B(K)$ 
by way of $U=\tilde\pi(z)$, and formula (\ref{pdEq1}) follows.  We omit the proof of 
the last sentence.  
\end{proof}

No operator theorist can resist repeating the elegant proof of 
von Neumann's inequality that flows from Theorem \ref{pdThm1}.  von Neumann's inequality 
\cite{vonNineq} 
asserts that for every operator $A\in\mathcal B(H)$ satisfying $\|A\|\leq 1$, one has 
\begin{equation}\label{pdEq2}
\|f(A)\|\leq \sup_{|z|\leq 1}|f(z)|
\end{equation}
for every polynomial $f(z)=a_0+a_1z+\cdots+a_nz^n$.  von Neumann's original proof was difficult, 
involving calculations with M{\"o}bius transformations and Blaschke products.  Letting 
$U\in\mathcal B(K)$ be a unitary power dilation of $A$ satisfying (\ref{pdEq1}), one has 
$f(A)=P_Hf(U)\restriction_H$, for every polynomial $f$, hence 
$$
\|f(A)\|\leq \|f(U)\|=\sup_{z\in\sigma(U)}|f(z)|\leq \sup_{|z|=1}|f(z)|.  
$$

\section{Completely positivity and Stinespring's theorem}\label{S:cp}

While one can argue that the GNS construction for states of \cstar s is a dilation theorem, 
it is probably best thought of as 
an application of the general method of associating a Hilbert space with a positive 
definite function as described in Section \ref{S:or}.  
Dilation theory {\em proper} went noncommutative in 1955 with the publication of a theorem 
of Stinespring \cite{stine}.  Stinespring once told me that his original 
motivation was simply to find a common generalization of 
Naimark's commutative result that a positive operator valued measure can be dilated to 
a spectral measure and the GNS construction for 
states of (noncommutative) \cstar s.  The theorem that emerged went well beyond that, and today 
has become a pillar upon which significant parts of operator 
theory and operator algebras rest.  The fundamental idea underlying 
the result was that of a completely positive linear map.  

The notion of positive linear functional or positive linear map is best thought of in  
a purely {\em algebraic} way.  More specifically, let $\mathcal A$ be a $*$-algebra, namely a 
complex algebra endowed with an antilinear mapping $a\mapsto a^*$ satisfying 
$(ab)^*=b^*a^*$ and $a^{**}=a$ for all $a,b\in \mathcal A$.  An
operator-valued linear map $\phi: \mathcal A\to \mathcal B(H)$ (and in particular a complex-valued 
linear functional $\phi: \mathcal A\to \mathbb C$) is called {\em positive} if it satisfies 
\begin{equation}\label{cpEq1}
\phi(a^*a)\geq 0, \qquad a\in \mathcal A. 
\end{equation}

One can promote this notion of positivity to matrix algebras over $\mathcal A$.  
For every $n=1,2,\dots$, the algebra $M_n(\mathcal A)$ of $n\times n$ matrices over $\mathcal A$ 
has a natural involution, in which the adjoint of an $n\times n$ matrix is defined as 
the transposed matrix of adjoints
$
(a_{ij})^*=(a_{ji}^*)
$, $1\leq i,j\leq n$, $a_{ij}\in \mathcal A$.  Fixing $n\geq 1$, 
a linear map $\phi: \mathcal A\to \mathcal B(H)$ 
induces a linear map $\phi_n$ from  $M_n(\mathcal A)$ to $n\times n$ operator matrices 
$(\phi(a_{ij}))$ which, after making the obvious identifications, can be viewed as a linear 
map of $M_n(\mathcal A)$ to operators on the direct sum of $n$ copies of $H$.  It makes 
good sense to say that $\phi_n$ is a positive linear map, and the original map 
$\phi$ is called {\em completely positive} if each $\phi_n$ is a positive linear map.  More 
explicitly, complete positivity at level $n$ requires that (\ref{cpEq1}) should hold 
for $n\times n$ matrices: For every 
$n\times n$ matrix $A=(a_{ij})$ with entries in $\mathcal A$ 
and every $n$-tuple of vectors $\xi_1,\dots, \xi_n\in H$, the $n\times n$ matrix 
$B=(b_{ij})$ defined by $B=A^*A$ satisfies 
$$
\sum_{i,j=1}^n\langle \phi(b_{ij})\xi_j,\xi_i\rangle 
=\sum_{i,j,k=1}^n\langle \phi(a_{ki}^*a_{kj})\xi_j,\xi_i\rangle \geq 0.  
$$
Note that this system of 
inequalities reduces to a somewhat simpler-looking system of inequalities (\ref{plEq1}) 
that we have already encountered in Section \ref{S:pl}.  

If $\mathcal A$ happens to be a \cstar,  then the elements $x\in \mathcal A$ that 
can be represented in the form $x=a^*a$ for some $a\in\mathcal A$ are precisely the 
self adjoint operators $x$ having nonnegative spectrum.  Since $M_n(\mathcal A)$ is also 
a \cstar\ in a unique way for every $n\geq 1$,  completely positive linear maps 
of \cstar s have a very useful spectral characterization: they should map self adjoint 
$n\times n$ operator matrices with nonnegative spectrum 
to self adjoint operators with nonnegative spectrum.    
Unfortunately, this spectral characterization 
breaks down completely for positive linear maps of more general Banach $*$-algebras, 
and in that more general 
context one must always refer back to positivity as it is 
expressed in (\ref{cpEq1}).  

Stinespring's original result was formulated in terms of operator maps 
defined on \cstar s.  We want to reformulate it somewhat 
into the more flexible context of linear maps of Banach $*$-algebras.  

\begin{thm}\label{cpThm1}
Let $\mathcal A$ be a Banach $*$-algebra with normalized unit and let 
$H$ be a Hilbert space.  For every 
completely positive linear map $\phi: \mathcal A\to\mathcal B(H)$ there is 
a representation $\pi: \mathcal A\to \mathcal B(K)$ of $\mathcal A$ on 
another Hilbert space $K$ and a bounded linear map $V: H\to K$ such that 
\begin{equation}\label{cpEq1.1}
\phi(a)=V^*\pi(a)V,\qquad a\in \mathcal A.  
\end{equation}
Moreover, the norm of the linking operator $V$ is given by $\|V\|^2=\|\phi(\mathbf 1)\|$.  
\end{thm}

We have omitted the statement and straightforward proof of the converse, namely 
that every linear map $\phi:\mathcal A\to \mathcal B(H)$ 
of the form (\ref{cpEq1.1}) must be completely positive, 
in order to properly emphasize the {\em construction} 
of the dilation from the basic properties of a completely positive map.

\begin{proof}  The underlying construction is identical with the original \cite{stine}, 
but a particular estimate requires care in the context of Banach $*$-algebras, 
and we will make that explicit.  
Consider the tensor product of complex vector spaces $\mathcal A\otimes H$, and 
the sesquilinear form $\langle \cdot,\cdot\rangle$ defined on it by setting 
$$
\langle \sum_{j=1}^ma_j\otimes\xi_j,\sum_{k=1}^n b_k\otimes\eta_k\rangle =
\sum_{j,k=1}^{m,n}\langle \phi(b_k^*a_j)\xi_j,\eta_k\rangle.  
$$
The fact that $\phi$ is completely positive implies that $\langle \zeta,\zeta\rangle\geq 0$ 
for every $\zeta\in \mathcal A\otimes H$.  
Letting $\mathcal N=\{\zeta\in \mathcal A\otimes H: \langle\zeta,\zeta\rangle=0\}$, 
the Schwarz inequality implies that $\mathcal N$ is a linear subspace and that 
the sesquilinear form can be promoted to an inner product on the quotient 
$K_0=(\mathcal A\otimes H)/\mathcal N$.  Let $K$ be the completion of the 
resulting inner product space.  

Each $a\in \mathcal A$ gives rise to a left multiplication operator $\pi(a)$ 
acting on $\mathcal A\otimes H$, defined uniquely by $\pi(a)(b\otimes \xi)=ab\otimes \xi$ 
for $b\in \mathcal A$ and $\xi\in H$.  The critical estimate that we require is
\begin{equation}\label{cpEq1.2}
\langle \pi(a)\zeta,\pi(a)\zeta\rangle\leq \|a\|^2\langle\zeta,\zeta\rangle, 
\qquad a\in \mathcal A, \quad \zeta\in \mathcal A\otimes H,  
\end{equation}
and it is proved as follows. Writing $\zeta=a_1\otimes\xi_1+\cdots+a_n\otimes\xi_n$, 
we find that 
\begin{align*}
\langle \pi(a)\zeta,\pi(a)\zeta\rangle&=
\sum_{j,k=1}^n\langle ab_j\otimes\xi_j,ab_k\otimes\xi_k\rangle
=\sum_{j,k=1}^n\langle\phi(b_k^*a^*ab_j)\xi_j,\xi_k\rangle\\
&=\sum_{j,k=1}^n\langle a^*ab_j\otimes\xi_j,b_k
\otimes \xi_k\rangle=\langle\pi(a^*a)\zeta,\zeta\rangle.  
\end{align*}
This formula implies that 
the linear functional $\rho(a)=\langle \pi(a)\zeta,\zeta\rangle$ 
satisfies $\rho(a^*a)=\langle \pi(a)\zeta,\pi(a)\zeta\rangle\geq 0$.   
Proposition \ref{a1Prop1} of the appendix implies 
$$
\rho(a^*a)\leq \rho(\mathbf 1)\|a^*a\|\leq\|\zeta\|^2\|a\|^2,   
$$ 
and (\ref{cpEq1.2}) follows.

It is obvious that $\pi(ab)=\pi(a)\pi(b)$ and that $\pi(\mathbf 1)$ is the identity 
operator.  Moreover, as in the argument above, we have 
$\langle \pi(a)\eta,\zeta\rangle=\langle \eta,\pi(a^*)\zeta\rangle$ for 
all $a\in \mathcal A$ and $\eta,\zeta\in \mathcal A\otimes H$.  
Finally, (\ref{cpEq1.2}) implies that $\pi(a)\mathcal N\subseteq\mathcal N$, 
so that each operator $\pi(a)$, $a\in\mathcal A$,  promotes naturally to a linear operator 
on the quotient $K_0=(\mathcal A\otimes H)/\mathcal N$.  Together with 
(\ref{cpEq1.2}), these formulas imply that $\pi$ gives rise 
to a $*$ representation of $\mathcal A$ as bounded operators on 
$K_0$ which extends uniquely 
to a representation of $\mathcal A$ on the completion 
$K$ of $K_0$, which we denote by the same letter $\pi$.  

It remains only to discuss the connecting operator $V$, which is defined 
by $V\xi=\mathbf 1\otimes \xi+\mathcal N$, $\xi\in H$.    One finds that 
$\pi(a)V\xi=a\otimes\xi+\mathcal N$, from which it follows that 
$$
\langle \pi(a)V\xi,V\eta\rangle=\langle a\otimes\xi+
\mathcal N,\mathbf 1\otimes\eta+\mathcal N\rangle
=\langle \phi(a)\xi,\eta\rangle, \qquad \xi,\eta\in H.  
$$
Taking $a=\mathbf 1$, we infer that $\|V\|^2=\|V^*V\|=\|\phi(\mathbf 1)\|$, and at that 
point the preceding formula implies $\phi(a)=V^*\pi(a)V$, $a\in \mathcal A$.  
\end{proof}

\section{operator spaces,  operator systems and extensions}\label{S:os}

In this section we discuss the basic features of operator spaces and their 
matrix hierarchies, giving only the briefest of overviews.  The 
interested reader is referred to one of the 
monographs \cite{BlLeMbook}, \cite{EffRu}, 
\cite{paulsenBk1} for more about this developing 
area of noncommutative analysis.  

Complex Banach spaces are the objects of a category whose maps 
are contractions -- linear operators of norm $\leq 1$.  The isomorphisms of this category are 
surjective isometries.  
A {\em function space} is a norm-closed linear subspace of some 
$C(X)$ -- the space of (complex-valued) continuous functions on a compact Hausdorff 
space $X$, endowed with the sup norm.  
All students of analysis know that every Banach space $E$ is isometrically 
isomorphic to a function space.  Indeed, the Hahn-Banach 
theorem implies that the natural map $\iota: E\to E^{\prime\prime}$ 
of $E$ into its double dual has 
the stated 
property after one views elements if $\iota(E)$ as continuous  
functions on the weak$^*$-compact unit ball $X$ of $E^\prime$.   
In this way the study of Banach spaces can be reduced to the study of 
function spaces, and that fact is occasionally useful.

An {\em operator space} is a norm-closed linear subspace $\mathcal E$ of 
the algebra $\mathcal B(H)$ of all bounded operators on a Hilbert 
space $H$.    Such an $\mathcal E$ is itself a Banach space, and is 
therefore isometrically isomorphic to a function space.    
However, the key fact about operator spaces is that they 
determine an entire hierarchy of operator spaces, one for every $n=1,2,\dots$.  
Indeed, for every $n$, the space $M_n(\mathcal E)$ of all $n\times n$ matrices over $\mathcal E$ 
is naturally an operator subspace of $\mathcal B(n\cdot H)$, $n\cdot H$ denoting 
the direct sum of $n$ copies of $H$.  Most significantly, a linear map 
of operator spaces $\phi:\mathcal E_1\to \mathcal E_2$ determines a sequence 
of linear maps $\phi_n: M_n(\mathcal E_1)\to M_n(\mathcal E_2)$, where $\phi_n$ 
is the linear map obtained by applying $\phi$ element-by-element to an 
$n\times n$ matrix over $\mathcal E_1$.  One says that $\phi$ is a 
{\em complete isometry} or a {\em complete contraction} if every 
$\phi_n$ is, respectively, an isometry or a contraction.  
There is a corresponding 
notion of {\em complete boundedness} that will not concern us here.  

Operator spaces can be viewed as the 
objects of a category whose maps are {\em complete} contractions.  The 
isomorphisms of this category are complete isometries, and one is led to 
seek properties of operator spaces that are invariant under this 
refined notion of isomorphism.  Like Shiva, a given Banach space acquires 
many inequivalent likenesses as an operator space.  And in operator space theory one 
pays attention to what happens at every level of the hierarchy.  The 
result is a significant and fundamentally noncommutative 
refinement of classical Banach space theory.

For example, since an operator space $\mathcal E\subseteq \mathcal B(H)$ 
is an ``ordinary" Banach space, it can be represented as a function 
system $\iota: \mathcal E\to C(X)$ as in the opening paragraphs of this 
section.  If we form the hierarchy of \cstar s $M_n(C(X))$, $n=1,2,\dots$, 
then we obtain a sequence of embeddings 
$$
\iota_n: M_n(\mathcal E)\to M_n(C(X)), \qquad n=1,2,\dots.  
$$
Note that the \cstar\ $M_n(C(X))$ is basically the \cstar\ of 
all matrix-valued continuous functions $F: X\to M_n(\mathbb C)$, with norm  
$$
\|F\|=\sup_{x\in X}\|F(x)\|, \qquad F\in M_n(C(X)).  
$$  
While the map $\iota$ is surely an isometry at the first level $n=1$, 
it may or may not be a complete isometry; indeed for the more 
interesting examples of operator spaces it is not.   Ultimately, the difference 
between these two categories can be traced to the noncommutativity 
of operator multiplication, and for that reason some analysts 
like to think of operator space theory as the ``quantized" reformulation of functional 
analysis.  

Finally, one can think of operator spaces somewhat more flexibly 
as norm-closed linear subspaces $\mathcal E$ 
of unital (or even nonunital) \cstar s $\mathcal A$.  That is 
because the hierarchy of \cstar s $M_n(\mathcal A)$ is well defined 
independently of any particular faithful realization 
$\mathcal A$ as a $C^*$-subalgebra of $\mathcal B(H)$.  

One can import the notion of {\em order} into the 
theory of operator spaces in a natural way.  A {\em function system} is 
a function space $E\subseteq C(X)$ with the property that 
$E$ is closed under complex conjugation and contains the 
constants.  One sometimes assumes that $E$ separates points 
of $X$ but we do not.  Correspondingly, an {\em operator system} 
is a self-adjoint operator space $\mathcal E\subseteq \mathcal B(H)$ 
that contains the identity operator $\mathbf 1$.  The natural 
notion of order between self adjoint operators, namely 
$A\leq B\iff B-A$ is a positive operator, has meaning in any 
operator system $\mathcal E$, and in fact every operator 
{\em system} is linearly spanned by its positive operators.  Every 
member $M_n(\mathcal E)$ of the matrix hierarchy over an operator system 
$\mathcal E$ is an operator system, so that it makes sense 
to speak of completely positive maps from one operator system 
to another.  

Krein's version of the Hahn-Banach theorem implies that  
a positive linear functional defined on an operator system 
$\mathcal E$ in a \cstar\ $\mathcal A$ can be extended to 
a positive linear functional on all of $\mathcal A$.  
It is significant that this extension theorem 
fails in general for operator-valued positive linear maps.  
Fortunately, the following result of 
\cite{arvSubalgI} provides an effective noncommutative 
counterpart of Krein's order-theoretic Hahn-Banach theorem: 

\begin{thm}\label{osThm1} Let $\mathcal E\subseteq \mathcal A$ be an operator system in a 
unital \cstar.  Then every 
operator-valued completely positive 
linear map $\phi: \mathcal E\to \mathcal B(H)$ can be extended 
to a completely positive linear map of $\mathcal A$ into $\mathcal B(H)$.  
\end{thm}
  
There is a variant of \ref{osThm1} that looks more like the 
original Hahn-Banach theorem.  Let $\mathcal E\subseteq \mathcal A$ 
be an operator space in a \cstar\ $\mathcal A$.  Then every 
operator-valued {\em complete} contraction $\phi:\mathcal E\to \mathcal B(H)$ 
can be extended to a completely contractive linear map of $\mathcal A$ to $\mathcal B(H)$.  
While the latter extension theorem emerged more than a decade after Theorem 
\ref{osThm1} (with a different and longer proof \cite{Wit1}, \cite{Witt2}), 
Vern Paulsen discovered a 
simple device that enables one to deduce it readily from the earlier result.  
That construction begins with an operator space $\mathcal E\subseteq \mathcal A$
and generates an associated operator {\em system} 
$\tilde{\mathcal E}$ in the $2\times 2$ matrix algebra $M_2(\mathcal A)$ over $\mathcal A$
as follows:   
$$
\tilde{\mathcal E}=\{
\begin{pmatrix}
\lambda\cdot\mathbf 1&A\\B^*&\lambda\cdot\mathbf 1
\end{pmatrix}
:A,B\in\mathcal E, \ \lambda\in \mathbb C\}.  
$$
Given a completely contractive linear map $\phi: \mathcal E\to \mathcal B(H)$, 
one can define a linear map $\tilde\phi: \tilde{\mathcal E}\to \mathcal B(H\oplus H)$ 
in a natural way
$$
\tilde\phi(
\begin{pmatrix}
\lambda\cdot\mathbf 1&A\\B^*&\lambda\cdot\mathbf 1
\end{pmatrix}
)=
\begin{pmatrix}
\lambda\cdot\mathbf 1&\phi(A)\\\phi(B)^*&\lambda\cdot\mathbf 1
\end{pmatrix}
, 
$$
and it is not hard to see that $\tilde\phi$ is completely positive (I have 
reformulated the construction in a minor but equivalent way for simplicity; see Lemma 
8.1 of \cite{paulsenBk2} for the original).  By Theorem 
\ref{osThm1}, $\tilde\phi$ extends to a completely positive 
linear map of $M_2(\mathcal A)$ to $\mathcal B(H\oplus H)$, and the behavior 
of that extension on the upper right corner is a  completely contractive 
extension of $\phi$.

\section{Spectral sets and higher dimensional operator theory}\label{S:ss}

Some aspects of commutative operator theory can be properly understood only when 
placed in the noncommutative context of the matrix hierarchies of the 
preceding section.  In this section we 
describe the phenomenon in concrete terms, referring the reader to the literature 
for technical details.  

Let $A\in \mathcal B(H)$ be a Hilbert space operator.  If $f$ is a rational function 
of a single complex variable that has no poles on the spectrum of $A$, then there is 
an obvious way to define an operator $f(A)\in \mathcal B(H)$.  Now fix a compact 
subset $X\subseteq \mathbb C$ 
of the plane that contains the spectrum of $A$.  The algebra $R(X)$ 
of all rational functions whose poles lie in the complement of $X$ forms a unital subalgebra 
of $C(X)$, and this functional calculus defines a unit-preserving homomorphism 
$f\mapsto f(A)$ of $R(X)$ into $\mathcal B(H)$.  One says that $X$ is a {\em spectral set} 
for $A$ if this homomorphism has norm $1$:
\begin{equation}\label{ssEq1}
\|f(A)\|\leq \sup_{z\in X}|f(z)|, \qquad f\in R(X).  
\end{equation}

von Neumann's inequality (\ref{pdEq2}) asserts that the closed unit disk is a spectral 
set for every contraction $A\in\mathcal B(H)$; indeed, that property is characteristic 
of contractions.  
While there is no corresponding characterization of the operators that have a more general 
set $X$ as a spectral set, we are still free to consider the 
class of operators that {\em do} have $X$ as 
a spectral set and ask if there is a generalization of Theorem \ref{pdThm1} that would 
apply to them.  Specifically, given an operator $A\in\mathcal B(H)$ that has $X$ as a spectral set, 
is there a normal operator $N$ acting on a larger Hilbert space $K\supseteq H$ 
such that the spectrum of $N$ is contained in the boundary $\partial X$ of $X$ and 
\begin{equation}\label{ssEq2}
f(A)=P_Hf(N)\restriction_H, \qquad f\in R(X)? 
\end{equation}

A result of Foias implies that the answer is yes if the complement of $X$ is connected, but 
it is no in general.  The reason the answer is no in general 
is that the hypothesis (\ref{ssEq1}) is not strong 
enough; and  that 
phenomenon originates in the {\em noncommutative} world. To see 
how the hypothesis must be strengthened, let $N$ be a normal operator 
with spectrum in $\partial X$.  The functional calculus for normal 
operators gives 
rise to a representation $\pi: C(\partial X)\to \mathcal B(K)$, 
$\pi(f)=f(N)$, $f\in C(\partial X)$.  It is easy to see that 
representations of \cstar s are completely positive and completely contractive linear 
maps, hence if the formula (\ref{ssEq2}) holds then the map $f\in R(X)\mapsto f(A)$ 
must be not only be a contraction, it must be a {\em complete} contraction.  

Let us examine the latter assertion 
in more detail.  Fix $n=1,2,\dots$ and let $M_n(R(X))$ be the algebra of all 
$n\times n$ matrices with entries in $R(X)$.  One can view an element of 
$M_n(R(X))$ as a matrix valued rational function 
$$
F: z\in X\mapsto F(z)=(f_{ij}(z))\in M_n(\mathbb C),
$$  
whose 
component functions belong to $R(X)$.  Notice that we can apply such a 
matrix valued function to an operator $A$ that has spectrum in $X$ to 
obtain an $n\times n$ matrix of operators -- or equivalently an 
operator $F(A)=(f_{ij}(A))$ in $\mathcal B(n\cdot H)$.  
The map $F\in M_n(R(X))\mapsto F(A)\in\mathcal B(n\cdot H)$ is a unit-preserving 
homomorphism of complex algebras.  $X$ is said to be a {\em complete} spectral 
set for an operator $A\in\mathcal B(H)$ if it contains the spectrum of $A$ 
and satisfies 
\begin{equation}\label{ssEq3}
\|F(A)\|\leq \sup_{z\in X}\|F(z)\|, \qquad F\in M_n(R(X)), \quad n=1,2,\dots.  
\end{equation}
Now if there is a normal operator $N$ with 
spectrum in $\partial X$ that relates to $A$ as in (\ref{ssEq2}), 
then for every $n=1,2,\dots$
$$
\|F(A)\|\leq \|F(N)\|\leq \sup_{z\in \partial X}\|F(z)\| = \sup_{z\in X}\|F(z)\|,  
$$
and we conclude that $X$ must be a {\em complete} spectral set for $A$.  

The following result from \cite{arvSubalgII} implies 
that complete spectral sets suffice for the existence of 
normal dilations.  It depends 
in an essential way on the extension theorem (Theorem \ref{osThm1}) for completely positive maps.   

\begin{thm}\label{ssThm1}
Let $A\in \mathcal B(H)$ be an operator 
that has a compact set $X\subseteq \mathbb C$ as a complete spectral set.  
Then there is a normal operator $N$ on 
a Hilbert space $K\supseteq H$ having spectrum in $\partial X$ such that 
$$
f(A)=P_Hf(N)\restriction_H, \qquad f\in R(X).  
$$
\end{thm}

The unitary power 
dilation of a contraction is unique up to natural equivalence.  That reflects 
a property of the unit circle $\mathbb T$: 
A positive linear map $\phi:C(\mathbb T)\to \mathcal B(H)$ 
is uniquely determined by its values on the nonnegative powers $1, z, z^2, \dots$ 
of the current variable $z$.   
In general, however, positive linear maps of $C(X)$ are not uniquely determined by 
their values on subalgebras of $C(X)$, with the result that there is 
no uniqueness assertion to complement the existence assertion of Theorem \ref{ssThm1} 
for the dilation theory of complete spectral sets.

 On the other hand, 
there is a ``many operators" generalization of Theorem \ref{ssThm1} that applies 
to completely contractive unit-preserving homomorphisms of arbitrary function algebras  
$A\subseteq C(X)$ that act on compact Hausdorff spaces $X$, in which $\partial X$ 
is replaced by the Silov boundary of $X$ relative to $A$.  The details can be found in Theorem 
1.2.2 of \cite{arvSubalgII} and its Corollary.

\section{Completely positive maps and endomorphisms}\label{S:en}

In recent years, certain problems arising in 
mathematical physics and quantum information theory have led researchers 
to seek a different kind of dilation theory, one that applies to 
semigroups of 
completely positive linear maps that act on von Neumann algebras.  In this section, we 
describe the simplest of these dilation theorems as it applies to the simplest semigroups 
acting on the simplest of von Neumann algebras.  A fuller accounting of these developments 
together with references to 
other sources can be found in Chapter 8 of the monograph \cite{arvMono}.    

Let $\phi: \mathcal B(H)\to \mathcal B(H)$ be a unit-preserving completely positive (UCP) map 
which is {\em normal} in the sense that for every normal state $\rho$ of $\mathcal B(H)$, 
the composition $\rho\circ\phi$ is also a normal state.  One can think of 
the semigroup 
$$
\{\phi^n: n=0,1,2,\dots\}
$$
as representing the discrete time evolution of an irreversible quantum system.  What we 
seek is another Hilbert space $K$ together with a normal 
$*$-endomorphism $\alpha: \mathcal B(K)\to \mathcal B(K)$ that is in some sense 
a ``power dilation" of $\phi$.  There are a number of ways one can make that vague idea 
precise, but only one of them is completely effective.  It is described as follows.

Let $K\supseteq H$ be a Hilbert space that contains $H$ and suppose we are 
given a normal $*$-endomorphism $\alpha: \mathcal B(K)\to\mathcal B(K)$ 
that satisfies $\alpha(\mathbf 1)=\mathbf 1$.  
We write the projection $P_H$ of $K$ on $H$ simply as $P$, and we  identify 
$\mathcal B(H)$ with the corner $P\mathcal B(K)P\subseteq \mathcal B(K)$.  
$\alpha$ is said to be a {\em dilation} of $\phi$ if 
\begin{equation}\label{enEq0}
\phi^n(A)=P\alpha^n(A)P, \qquad A\in \mathcal B(H)=P\mathcal B(K)P, \quad n=0,1,2,\dots.  
\end{equation}

Since $\phi$ is a unit-preserving map of $\mathcal B(H)$,  $P=\phi(P)=P\alpha(P)P$, 
so that $\alpha(P)\geq P$.  Hence 
we obtain an increasing sequence of projections 
\begin{equation}\label{enEq1}
P\leq \alpha(P)\leq \alpha^2(P)\leq\cdots.  
\end{equation}
The  limit projection 
$P_\infty=\lim_n \alpha^n(P)$ satisfies $\alpha(P_\infty)=P_\infty$, hence the compression of 
$\alpha$ to the larger corner $P_\infty\mathcal B(K)P_\infty\cong\mathcal B(P_\infty K)$ 
of $\mathcal B(K)$ is a unital 
$*$-endomorphism that is itself a dilation of $\phi$.  By cutting down if necessary 
we can assume that the configuration is {\em proper} in the sense that 
\begin{equation}\label{enEq2}
\lim_{n\to \infty}\alpha^n(P)=\mathbf 1_K,  
\end{equation}
and in that case the endomorphism $\alpha$ is said to be a {\em proper} dilation 
of $\phi$.  We have refrained from using the term {\em minimal} to describe this 
situation because in the context of semigroups of completely positive maps, 
the notion of {\em minimal} dilation is a more subtle one that requires 
a stronger hypothesis.  That hypothesis is discussed in Remark \ref{enRem2} below.  

\begin{rem}[Stinespring's theorem is not enough]\label{enRem1}  It is by no means obvious 
that dilations should exist.  One might attempt to construct a dilation of the semigroup 
generated by a single UCP map 
$\phi: \mathcal B(H)\to \mathcal B(H)$ by applying Stinespring's theorem to the individual 
terms of the 
sequence of powers $\phi^n$, $n=0,1,2,\dots$, and then somehow putting the pieces together 
to obtain the dilating endomorphism $\alpha$.  Indeed, Stinespring's 
theorem provides us with a Hilbert space $K_n\supseteq H$ and a 
representation $\pi_n: \mathcal B(H)\to \mathcal B(K_n)$ for every $n\geq 0$ such that 
$$
\phi^n(A)=P_H\pi_n(A)\restriction_H, \qquad A\in\mathcal B(H), \quad n=0,1,2,\dots.  
$$
However, while these formulas certainly inherit a relation to each other by virtue of 
the semigroup formula 
$\phi^{m+n}=\phi^m\circ\phi^n$, $m,n\geq 0$, if one attempts to 
exploit these relationships one finds 
that the relation between 
$\pi_m$, $\pi_n$ and $\pi_{m+n}$ is extremely 
awkward.  Actually, 
there is no apparent way to assemble the von Neumann algebras $\pi_n(\mathcal B(H))$ 
into a single von Neumann algebra that plays the role of $\mathcal B(K)$, 
on which one can define a single endomorphism $\alpha$ that converts these formulas 
into the single formula (\ref{enEq0}).  
Briefly put, 
{\em Stinespring's theorem does not apply to semigroups}.  
\end{rem}

These observations suggest that 
the problem of constructing dilations in this context should require an entirely new method, 
and it does.  
The proper result for normal UCP maps acting on 
$\mathcal B(H)$ was discovered by Bhat and Parthasarathy \cite{BhatPar}, building 
on earlier work of Parthasarathy \cite{parBookQProb} that was set 
in the context of quantum probability theory.  The result was later 
extended by Bhat to semigroups of completely positive maps that act 
on arbitrary von Neumann algebras \cite{bhatMin}.  The construction of the dilation 
has been reformulated in various 
ways; the one I like is in Chapter 8 of \cite{arvMono} (also see \cite{arvNCgen}).  
Another approach, due to Muhly 
and Solel \cite{muhSol2}, is based on correspondences over von Neumann algebras.  The history 
of earlier approaches to 
this kind of dilation theory is summarized in the notes of Chapter 8 of \cite{arvMono}.

We now state the appropriate result for $\mathcal B(H)$ without proof:  

\begin{thm}\label{enThm1}
For every normal UCP map $\phi:\mathcal B(H)\to \mathcal B(H)$, there is a Hilbert 
space $K\supseteq H$ and a normal 
$*$-endomorphism $\alpha:\mathcal B(H)\to \mathcal B(H)$ satisfying $\alpha(\mathbf 1)=\mathbf 1$ 
that is a proper dilation 
of $\phi$ as in (\ref{enEq0}).  
\end{thm}

\begin{rem}[Minimality and uniqueness]\label{enRem2}  The notion of minimality for a dilation 
$\alpha:\mathcal B(K)\to \mathcal B(K)$ of a UCP map 
$\phi: \mathcal B(H)\to \mathcal B(H)$ is described as follows.  Again, we identify 
$\mathcal B(H)$ with the corner $P\mathcal B(K)P$.  We have already pointed out that 
the projections $\alpha^n(P)$ increase with $n$.  However, the sequence of (nonunital)
von Neumann subalgebras $\alpha^n(\mathcal B(H))$, $n=0,1,2,\dots$,  neither increases 
nor decreases with $n$, and 
that behavior requires care.   The proper notion of minimality in this context is that the 
set of all 
vectors in $K$ of the form 
\begin{equation}\label{enEq4}
\alpha^{n_1}(A_1)\alpha^{n_2}(A_2)\cdots \alpha^{n_k}(A_k)\xi, 
\end{equation}
where $k=1,2,\dots$, $n_k=0,1,2,\dots$, $A_k\in \mathcal B(H)$, and $\xi\in H$, should have 
$K$ as their closed linear span.  Equivalently, 
the smallest subspace of $K$ that contains $H$ and is invariant under the set of operators 
$$
\mathcal B(H)\cup \alpha(\mathcal B(H))\cup\alpha^2(\mathcal B(H))\cup\cdots
$$ 
should be all of $K$.  It is a fact that every minimal dilation is proper, but the 
converse is false.  It is also true that every proper dilation can be reduced in 
a natural way to a minimal one, and finally, that any two minimal dilations of the 
semigroup $\{\phi^n: n\geq 0\}$  are isomorphic in a natural sense.  

We also point out that there is 
a corresponding dilation theory for one-parameter semigroups of UCP maps.  These facts are 
discussed at length in Chapter 8 of \cite{arvMono}.  
\end{rem}

\appendix
\section{Brief on Banach $*$-algebras}\label{S:a1}

Banach $*$-algebras (defined at the beginning of Section \ref{S:pl})
are useful because they are flexible -- it is usually 
a simple matter to 
define a Banach $*$-algebra with the properties one needs.  More importantly, 
it is far easier to define states and representations of Banach 
$*$-algebras than it is for the more rigid category of \cstar s.  
For example, we made use of the technique in the proof of Theorem \ref{apThm1}
and the estimate of Remark \ref{apRem3}.  

On the other hand, it is obviously desirable to have $C^*$-algebraic tools 
available for carrying out analysis.  
Fortunately one can have it both ways, because every Banach $*$-algebra $\mathcal A$ 
is associated with a unique enveloping \cstar\ $C^*(\mathcal A)$ 
which has the ``same" representation theory 
and the ``same" state space as $\mathcal A$.  In this Appendix we briefly describe 
the properties of 
this useful functor $\mathcal A\to C^*(\mathcal A)$ for the 
category of Banach $*$-algebras that have a normalized unit $\mathbf 1$.  
There are similar results (including Proposition \ref{a1Prop1} below) 
for many nonunital Banach 
$*$-algebras -- including the group algebras of locally compact groups --  
provided that they have appropriate approximate units.  A comprehensive 
treatment can be found in  \cite{dixCstar}.

The fundamental fact on which these results are based is the 
following (see Proposition 4.7.1 of the text \cite{arvSpecTh} for a proof):  

\begin{prop}\label{a1Prop1}  Every positive linear functional $\rho$ on a unital 
Banach $*$-algebra $\mathcal A$ is bounded, and in fact 
$\|\rho\|=\rho(\mathbf 1)$.  
\end{prop}

What we actually use here is the following consequence, which is proved by applying 
Proposition \ref{a1Prop1} to functionals of the form $\rho(a)=\langle\phi(a)\xi,\xi\rangle$:    

\begin{cor}\label{a1Cor1}
Every operator-valued positive linear 
map $\phi: \mathcal A\to \mathcal B(H)$ is bounded, and $\|\phi\|=\|\phi(\mathbf 1)\|$.  
\end{cor}

By a {\em representation} of a Banach $*$-algebra $\mathcal A$ we mean 
a $*$-preserving homomorphism $\pi: \mathcal A\to \mathcal B(H)$ of 
$\mathcal A$ into the $*$-algebra of operators on a Hilbert space.  It is useful 
to assume the representation is nondegenerate in the 
sense that $\pi(\mathbf 1)=\mathbf 1$; if that is not the case, it can 
be arranged by passing to the subrepresentation defined on 
the subspace $H_0=\pi(\mathbf 1)H$.  Representations of Banach $*$-algebras 
 arise from positive linear functionals 
(by way of the GNS construction which makes use of Proposition \ref{a1Prop1}) or from 
completely positive linear maps (by a variation of Theorem \ref{cpThm1}, by making use 
of Corollary \ref{a1Cor1}).

While we have made no hypothesis on the norms $\|\pi(a)\|$ associated with 
a representation $\pi$, it follows immediately from Proposition \ref{a1Prop1} that every 
representation of $\mathcal A$ has norm $1$.  Indeed, for every 
unit vector $\xi\in H$ and $a\in \mathcal A$, 
$\rho(a)=\langle \pi(a)\xi,\xi\rangle$ defines a positive 
linear functional on $\mathcal A$ with $\rho(\mathbf 1)=1$, so that 
$$
\|\pi(a)\xi\|^2=\langle\pi(a)^*\pi(a)\xi,\xi\rangle =\langle\pi(a^*a)\xi,\xi\rangle
=\rho(a^*a)\leq \|a^*a\|\leq  \|a\|^2,   
$$ 
and  $\|\pi(a)\|\leq \|a\|$ follows.   It is an instructive exercise 
to find a direct proof of the 
inequality $\|\pi(a)\|\leq \|a\|$ that does not make use of Proposition \ref{a1Prop1}.  

\begin{rem}[Enveloping $C^*$-algebra of a Banach $*$-algebra]\label{a1Rem1}
Consider the seminorm $\|\cdot\|_1$ defined on $\mathcal A$ by 
$$
\|a\|_1=\sup_\pi\|\pi(a)\|, \qquad a\in \mathcal A, 
$$
the supremum taken over a ``all" representations of $\mathcal A$.
Since the representations of $\mathcal A$ do not form a set, the 
quotes simply refer to an obvious way of choosing sufficiently 
many representatives from unitary equivalence classes 
of representations so that 
every representation is unitarily equivalent to a direct sum 
of the representative ones.  It is clear that 
$\|a^*a\|_1=\|a\|_1^2$.  Indeed, $\|\cdot\|_1$ is a $C^*$-seminorm, 
and the completion of $\mathcal A/\{x\in \mathcal A: \|x\|_1=0\}$
is a \cstar\ $C^*(\mathcal A)$, called the enveloping \cstar\ of 
$\mathcal A$.  The natural completion map 
\begin{equation}\label{a1Eq1}
\iota:  \mathcal A\to C^*(\mathcal A)
\end{equation}
is a $*$-homomorphism having dense range and norm $1$.  
This completion 
(\ref{a1Eq1}) has the following 
universal property: For every representation $\pi: \mathcal A\to \mathcal B(H)$ 
there is a unique representation $\tilde\pi: C^*(\mathcal A)\to \mathcal B(H)$ 
such that $\tilde\pi\circ\iota=\pi$.  The map 
$\pi\to \tilde\pi$ is in fact a bijection.  Indeed,  Proposition \ref{a1Prop1} 
is equivalent to the assertion that there is a bijection 
between the set of positive linear functionals $\rho$ 
on $\mathcal A$ and the set of positive linear functionals $\tilde\rho$ 
on its enveloping \cstar,  defined by a similar formula 
$\tilde\rho\circ\iota=\rho$.

One should keep in mind that the completion map (\ref{a1Eq1}) can have a 
nontrivial kernel in general, but for many important examples it is injective.  
For example, 
it is injective in the case of 
group algebras -- the Banach $*$-algebras $L^1(G)$ associated with locally compact 
groups $G$.  When $G$ is commutative, the enveloping \cstar\ of 
$L^1(G)$ is the \cstar\ $C_\infty(\hat G)$ 
of continuous functions that vanish at $\infty$ on the character group 
$\hat G$ of $G$.  
\end{rem}

\bibliographystyle{alpha}

\newcommand{\noopsort}[1]{} \newcommand{\printfirst}[2]{#1}
  \newcommand{\singleletter}[1]{#1} \newcommand{\switchargs}[2]{#2#1}

\end{document}